\documentclass[12pt]{amsart}
\usepackage{amscd}
\def\frk{\frak}               

\def\Phi{{\frk n}}
\def\Phi{{\frk N}}
\def\opn#1#2{\def#1{\operatorname{#2}}} 
\opn\chara{char} \opn\length{\ell} \opn\pd{pd} \opn\rk{rk}
\opn\projdim{proj\,dim} \opn\injdim{inj\,dim} \opn\rank{rank}
\opn\depth{depth} \opn\grade{grade} \opn\height{height}
\opn\embdim{emb\,dim} \opn\codim{codim}

\opn\Tr{Tr} \opn\bigrank{big\,rank}
\opn\superheight{superheight}\opn\lcm{lcm}
\opn\trdeg{tr\,deg}
\opn\reg{reg} \opn\lreg{lreg} \opn\ini{in} \opn\lpd{lpd}
\opn\size{size}
\opn\div{div} \opn\Div{Div} \opn\cl{cl} \opn\Cl{Cl}
\opn\Spec{Spec} \opn\Supp{Supp} \opn\supp{supp} \opn\Sing{Sing}
\opn\Ass{Ass} \opn\Min{Min}
\opn\Ann{Ann} \opn\Rad{Rad} \opn\Soc{Soc}
\opn\Im{Im} \opn\Ker{Ker} \opn\Coker{Coker} \opn\Am{Am}
\opn\Hom{Hom} \opn\Tor{Tor} \opn\Ext{Ext} \opn\End{End}
\opn\Aut{Aut} \opn\id{id}

\opn\nat{nat}
\opn\pff{pf}
\opn\Pf{Pf} \opn\GL{GL} \opn\SL{SL} \opn\mod{mod} \opn\ord{ord}
\opn\Gin{Gin} \opn\Hilb{Hilb}
\opn\aff{aff} \opn\con{conv} \opn\relint{relint} \opn\st{st}
\opn\lk{lk} \opn\cn{cn} \opn\core{core} \opn\vol{vol}
\opn\link{link} \opn\star{star}
\opn\gr{gr}

\def\pot#1#2{#1[\kern-0.28ex[#2]\kern-0.28ex]}

\opn\dirlim{\underrightarrow{\lim}}
\opn\inivlim{\underleftarrow{\lim}}
\def\Implies{\ifmmode\Longrightarrow \else
        \unskip${}\Longrightarrow{}$\ignorespaces\fi}
\def\implies{\ifmmode\Rightarrow \else
        \unskip${}\Rightarrow{}$\ignorespaces\fi}
\def\iff{\ifmmode\Longleftrightarrow \else
        \unskip${}\Longleftrightarrow{}$\ignorespaces\fi}

\let\:=\colon
\newtheorem{Theorem}{Theorem}[section]

\let\epsilon\varepsilon
\let\phi=\varphi
\let\kappa=\varkappa
\def\qed{\ifhmode\textqed\fi
      \ifmmode\ifinner\quad\qedsymbol\else\dispqed\fi\fi}
\def\textqed{\unskip\nobreak\penalty50
       \hskip2em\hbox{}\nobreak\hfil\qedsymbol
       \parfillskip=0pt \finalhyphendemerits=0}
\def\dispqed{\rlap{\qquad\qedsymbol}}

\opn\dis{dis}
\def\pnt{{\raise0.5mm\hbox{\large\bf.}}}

\opn\Lex{Lex}

\begin{document}

\title{On the Ramsey numbers for a combination of paths and Jahangirs}

\author{Kashif Ali, Edy Tri Baskoro}

\address{Kashif Ali, School of Mathematical Sciences, 68-B New Muslim Town,
    Lahore,
    Pakistan.}\email{akashifali@gmail.com}
\address{Edy Tri Baskoro, Combinatorial Mathematics Research Group,
Faculty of Mathematics and Natural Sciences, Institut Teknologi
Bandung, Jalan Genesa 10 Bandung 40132, Indonesia.}
\email{ebaskoro@math.itb.ac.id} \maketitle

\begin{abstract}
For given graphs $G$ and $H,$ the \emph{Ramsey number} $R(G,H)$ is
the least natural number $n$ such that for every graph $F$ of order
$n$ the following condition holds: either $F$ contains $G$ or the
complement of $F$ contains $H.$ In this paper, we improve the
Surahmat and Tomescu's result \cite{ST:06} on the Ramsey number of
paths versus Jahangirs. We also determine the Ramsey number $R(\cup
G,H)$, where $G$ is a path and $H$ is a Jahangir graph.
\end{abstract}

\section{Introduction}
The study of Ramsey Numbers for (general) graphs have received
tremendous efforts in the last two decades, see few related papers
\cite{BBS:05,Bur:75,BRR:88,CZZ:04,HasB:04,SB:01} and a nice survey
paper \cite{Rad:04}. One of useful results on this is the
establishment of a general lower bound by Chv\'atal and Harary
\cite{Catal:72}, namely $R(G,H)\geq (\chi (G)-1)(c(H)-1)+1$, where
$\chi (G)$ is the chromatic number of $G$ and $c(H)$ is the number
of vertices in
the largest component of $H$.\\

\noindent Let $G(V,E)$ be a graph with the vertex-set $V(G$) and
edge-set $E(G)$. If $(x,y) \in E(G)$ then $x$ is called {\em
adjacent} to $y$, and $y$ is a {\em neighbor} of $x$ and vice versa.
For any $A \subseteq V(G)$, we use $N_A(x)$ to denote the set of all
neighbors of $x$ in $A$, namely $N_A(x) = \{ y \in A | (x,y) \in
E(G)\}$. Let $P_{n}$ be a path with $n$ vertices, $C_n$ be a cycle
with $n$ vertices, and $W_m$ be a wheel of $m+1$ vertices, i.e., a
graph consisting of a cycle $C_m$ with one additional vertex
adjacent to all vertices of $C_m$. For $m \geq 2$, the {\em Jahangir
graph} $J_{2m}$ is a graph consisting of a cycle $C_{2m}$ with one
additional vertex adjacent alternatively to $m$ vertices of
$C_{2m}$. For example, Figure 1\footnote{The figure $J_{16}$ appears
on Jahangir's tomb in his mausoleum, it lies in 5 km north-west of
Lahore, Pakistan across the River Ravi. His tomb was built by his
Queen Noor Jehan and his son Shah-Jehan (This was emperor who
constructed one of the wonder of world Taj Mahal in India) around
1637 A.D. It has a majestic structure made of red sand-stone and
marble.} \\

\noindent Recently, Surahmat and Tomescu \cite{ST:06} studied the
Ramsey number
of a combination of $P_n$ versus a $J_{2m}$, and obtained the following result.\\

\noindent {\bf Theorem~A}. \cite{ST:06}.\\
 {\it $R(P_n,J_{2m})=\left\{
\begin{array}{ll}
6&\mbox{if $(n,m)=(4,2)$,}\\
n+1& \mbox{if $m=2$ and $n\geq 5$,}\\
 n+m-1& \mbox{~if $m\geq 3$ and $n\geq (4m-1)(m-1)+1$.}
\end{array}\right.$}\\

\noindent In this paper, we determine the Ramsey numbers involving
paths and Jahangir graphs. For particular, we improve the Surahmat
and
Tomescu's result for Jahangir graphs $J_6$, $J_8$ and  $J_{10}$ as follows.\\

\begin{Theorem}
$R(P_n,J_{2m})=n+m-1$ for $ n \geq 2m+1$ and $m=3,4$ or $5$.
\end{Theorem}

\noindent We are also able to determine the Ramsey number $R(kP_n,
J_{2m})$, for any integer $k \geq 2$, $m \geq 2$.
 These results are stated in the following theorems.

\begin{Theorem}
 $R(kP_n,J_4)=kn+1$, for $n \geq 4$, $k \geq 1$, except for $(n=4,k=1)$.
\end{Theorem}

\begin{Theorem}
 $R(kP_n,J_{2m})=kn+m-1$, for any integer $n \geq 2m+1$ if $m=3,4$ or $5$;
 and for $n\geq(4m-1)(m-1)+1$ if $m \geq 6$, where $k \geq 2$.
\end{Theorem}

\section{The Proof of Theorems}

\noindent {\bf The proof of Theorem 1}.\\
Consider graph $G \cong K_{m-1} \cup K_{n-1}$. Clearly, $G$ contains
no $P_n$ and $\overline{G}$ contains no $J_{2m}$. Thus,
$R(P_n,J_{2m})\geq n+m-1$. For $m=3,4$ or $5$ and $n \geq 2m+1$, we
will show that $R(P_n,J_{2m})\leq n+m-1$. Let $F$ be a graph of
$n+m-1$ vertices containing no $P_n$. Take any longest path $L$ in
$F$. Let $L$ be  $(x_1, x_2, \cdots, x_k)$, and $Y=V(F)\backslash
V(L)$. Since $k \leq n-1$, then $|Y| \geq m$. Obviously, $yx_1,yx_k$
are not in $E(F)$, for any $y \in Y$. Now, consider the following two cases\\

\noindent{\bf Case 1.} $2m \leq |L|\leq n-1$. \\
Let $|L|=t$ and $A=\{x_2,x_3, \cdots ,x_{2m-1}\}$ be the set of
first $2m-2$ vertices of $L$ after $x_1$. Take the set of any $m$
distinct vertices of $Y$ and denote it by $B=\{y_1,\cdots,y_m\}$. By
the maximality of $L$, every vertex of $B$ has at most $m-1$
neighbors in $A$. If there are two vertices of $B$ having $m-1$
neighbors in $A$ then all the
neighbors are intersected.\\

\noindent {\bf Subcase 1.1} There exists $b \in B$, $|N_A(b)|= m-1$.\\
Let $A_1=A\backslash N_A(b)$ and take any vertex $v_1$ of $A_1$
whose the highest degree at $B$. Define $D_1=\{x_1,x_t,b\} \cup
A_1\backslash \{v_1\}$, and $D_2=\{v_1\} \cup B\backslash\{b\}$. By
the maximality of $L$, $d_{D_1}(w) \leq 1$ for any vertex $w$ of
$D_2$. In particular, $d_{D_1}(v_1)=0$. Since $v_1$ has the highest
degree then there are at most $m-2$ edges connecting vertices
between
$D_1$ and $D_2$ in $F$. This implies that $D_1 \cup D_2$ will induces a $J_{2m}$ in $\overline{F}$.\\

\noindent {\bf Subcase 1.2} All vertices $b \in B$, $|N_A(b)|\leq m-2$.\\
If $m=3$ then let $D_1 =\{$any two vertices of $A\}$. If $m=4$ then
by the pigeonhole principle there exists two vertices of $A$ has
neighbors at most $1$ in $B$. In this case let $D_1=\{$three
vertices of $A$ with two of degree at most one $\}$. If $m=5$ then
by the pigeonhole principle there exists three vertices of $A$ has
neighbors at most $2$ in $B$. In this case let $D_1=\{$four vertices
of $A$ with three of degree at most two$\}$.
Therefore, $\{x_1,x_t\} \cup D_1 \cup B$ will induce a $J_{2m}$ in $\overline{F}$.\\

\noindent{\bf Case 2.} $1\leq |L|\leq 2m-1$.\\
We breakdown the proof into several subcases.\\

\noindent {\bf Subcase 2.1}. $1 \leq |L| \leq 3$\\
In this case, the component of $F$ is either $K_1$, $P_2$,
$C_3$ or a star. Therefore, $\overline{F}$ contains a $J_{2m}$, for $m=3,4$ or 5.\\

\noindent {\bf Subcase 2.2}. $4 \leq |L| \leq m+1$.\\
Let $L$ be $(x_1,x_2, \cdots, x_{t})$, where $t \leq m+1$, and so
$|Y|=|V(F)\backslash V(L)|$ $\geq 2m-1$. Now, consider the set
$N_Y(x_2)$ of vertices in $Y$ adjacent to $x_2$. Note that any
vertex of $N_Y(x_2)$ is nonadjacent to any other vertices of $Y$. If
$|N_Y(x_2)| \geq m-2$ then form two sets $D_1$ and $D_2$ as follows.
The set $D_1$ consists of $x_1, x_t$ and any $m-2$ vertices of
$N_Y(x_2)$. The set $D_2$ consists of the other vertices of $Y$ not
selected in $D_1$. Thus, $|D_1|=m$ and $|D_2|=m+1$. By the
maximality of $L$, there is no edge connecting any vertex of $D_1$
to any vertex of $D_2$. Thus, the set $D_1 \cup D_2$ induces
$K_{m,m+1} \supseteq J_{2m}$ in $\overline{F}$. If $|N_Y(x_2)|= m-3$
then  take $D_1=\{x_1,x_t,x_2\} \cup N_Y(x_2)$, and $D_2$ as the set
of the remaining vertices of $Y$. Then, $D_1 \cup D_2$ again
contains $K_{m,m+1} \supseteq J_{2m}$ in $\overline{F}$. Now, if
$|N_Y(x_2)|= m-4$ (for $m=4$ or $5$) then in showing $\overline{F}
\supseteq J_{2m}$ take $D_1=\{x_1,x_t,x_2,x_{t-1}\} \cup N_Y(x_2)$,
and $D_2$ as the set of the remaining vertices of $Y$ not adjacent
to $x_{t-1}$. This is true since $|N_Y(x_{t-1})| \leq 1$ (by
symmetrical argument). If $|N_Y(x_2)|= m-5$ (for $m=5$ only), then
$D_1=\{x_1,x_2,x_{t-1},x_t,b\}$ where $b$ is a vertex at distance
two from $x_3$ or $b$ is any vertex of $Y$ with a smallest degree,
and $D_2$ as the set
of the remaining vertices of $Y$. Thus, $D_1 \cup D_2$ will induce $J_{10}$ in $\overline{F}$.\\

\noindent {\bf Subcase 2.3}. $|L|=m+2$.\\
Let $L$ be $(x_1,x_2,\cdots, x_{t})$ where $t=m+2$, then
$|Y|=|V(F)\backslash V(L)| \geq 2m-2$. Now, consider the set
$N_Y(x_2)$ of vertices in $Y$ adjacent to $x_2$. Note that any
vertex of $N_Y(x_2)$ is nonadjacent to any other vertices of $Y$. If
$|N_Y(x_2)| \geq m-2$ then form two sets $D_1$ and $D_2$ as follows.
If $x_3$ is nonadjacent to $x_{m+2}$ then $D_1=\{x_1,x_{m+2}\} \cup
\{$any $m-2$ vertices of $N_Y(x_2)\}$ and $D_2$ consists of $x_3$
together with the remaining vertices of $Y$. Otherwise (if $x_3 \sim
x_{m+2}$), take $D_1=\{x_1,x_{m+2},x_4\} \cup \{$any $m-2$ vertices
of $N_Y(x_2)\}$ and $D_2$ consists of any $m$ remaining vertices of
$Y$. By the maximality of $L$, there is no edge connecting any
vertex of $D_1$ to any vertex of $D_2$.
Thus, the set $D_1 \cup D_2$ induces $K_{m,m+1} \supseteq J_{2m}$ in $\overline{F}$.\\

\noindent If $|N_Y(x_2)|= m-3$ then  take $D_1=\{x_1,x_t,x_2\} \cup
N_Y(x_2)$, and $D_2$ as the set of the remaining vertices of $Y$.
Then, $D_1 \cup D_2$ again contains $K_{m,m+1} \supseteq J_{2m}$ in
$\overline{F}$. Now, if $|N_Y(x_2)|= m-4$ (for $m=4$ or $5$) then in
showing $\overline{F} \supseteq J_{2m}$ take
$D_1=\{x_1,x_t,x_2,x_{t-1}\} \cup N_Y(x_2)$, and $D_2$ as the set of
the remaining vertices of $Y$ not adjacent to $x_{t-1}$. This is
true since $|N_Y(x_{t-1})| \leq 1$ (by symmetrical argument). If
$|N_Y(x_2)|= m-5$ (for $m=5$ only), then
$D_1=\{x_1,x_2,x_{t-1},x_t,b\}$ where $b$ is a vertex at distance
two from $x_3$ or $b$ is any vertex of $Y$ with a smallest degree,
and $D_2$ as the set
of the remaining vertices of $Y$. Thus, $D_1 \cup D_2$ will induce $J_{10}$ in $\overline{F}$.\\

\noindent {\bf Subcase 2.4}. $|L|=m+3$ (or $2m-1, 2m-2$ if $m=4,5$ respectively).\\
Let $L$ be $(x_1,x_2,\cdots, x_{t})$ where $t=m+3$, then
$|Y|=|V(F)\backslash V(L)| \geq 2m-3$. Now, consider the set
$N_Y(x_2)$ of vertices in $Y$ adjacent to $x_2$. Note that any
vertex of $N_Y(x_2)$ is nonadjacent to any other vertices of $Y$. If
$|N_Y(x_2)| \geq m-1$ then form two sets $D_1$ and $D_2$ as follows.
If $x_{t-1}$ is adjacent to some vertex of $N_Y(x_2)$ then by the
maximality of $L$, $x_{t-2}$ is nonadjacent to $x_1$ and any vertex
of $N_Y(x_2)$. In this case set $b=x_{t-2}$.  If $x_{t-1}$ is
nonadjacent to any vertex of $N_Y(x_2)$, then take $b=x_{t-1}$
provided $x_{t-1} \not\sim x_1$. Otherwise (if $x_{t-1}\sim x_1$),
by the maximality of $L$ we have that $x_{t-2}$ is nonadjacent to
$x_1$ and to any vertex of $N_Y(x_2)$. In this case, again take
$b=x_{t-2}$. Now, define $D_1=\{x_1\} \cup \{$any $m-1$ vertices of
$N_Y(x_2)\}$ and $D_2=\{x_3,x_t, b\} \cup$  $\{$ any  $m-2$ other
vertices of $Y\}$. By the maximality of $L$, there is no edge
connecting any vertex of $D_1$ to any vertex of $D_2$.
Thus, the set $D_1 \cup D_2$ induces $K_{m,m+1} \supseteq J_{2m}$ in $\overline{F}$.\\

\noindent If $|N_Y(x_2)|= m-2$ then  take $D_1=\{x_1,x_2\} \cup
N_Y(x_2)$, and $D_2 =\{x_3,x_t\} \cup \{$ any $m-1$ other vertices
of $Y$. Then, $D_1 \cup D_2$ contains $K_{m,m+1}$ minus at most two
edges $(x_2,x_3)$ and $(x_2,x_t)$ in $\overline{F}$. Therefore,
$\overline{F} \supseteq J_{2m}$. Now, if $|N_Y(x_2)|= m-3$ then in
showing $\overline{F} \supseteq J_{2m}$ take $D_1=\{x_1,x_2,x_t\}
\cup N_Y(x_2)$, and $D_2 =\{x_3\} \cup \{$ any $m$ other vertices of
$Y$. This is true since $D_1 \cup D_2$ contains $K_{m,m+1}$ minus at
most two edges $(x_2,x_3)$ and $(x_2,x_t)$ in $\overline{F}$. If
$|N_Y(x_2)|= m-4$, then $D_1=\{x_1,x_2,x_{t-1},x_t\} \cup N_Y(x_2)
\cup N_Y(x_{t-1})$ and $D_2$ as the set of the remaining vertices of
$Y$. Thus, $D_1 \cup D_2$ will induce $K_{m,m+1}$ in $\overline{F}$.
If $|N_Y(x_2)|= m-5$ (only for $m=5$), then
$D_1=\{x_1,x_2,x_{t-1},x_t,b\}$, where $b$ is either $x_3$, a
neighbor of $x_3$ in $Y$ or a vertex of $Y$ at distance two from
$x_3$  and $D_2$ as the set
of the remaining vertices of $Y$. Thus, $D_1 \cup D_2$ will induce $K_{m,m+1}$ minus at most one edge in $\overline{F}$.\\


\noindent {\bf Subcase 2.5}. $|L|=m+4=2m-1$ (only for $m=5$).\\
Let $L$ be $(x_1,x_2,\cdots, x_{t})$ where $t=2m-1$, then
$|Y|=|V(F)\backslash V(L)| \geq 2m-4$. Now, consider the set
$N_Y(x_2)$ of vertices in $Y$ adjacent to $x_2$. Note that any
vertex of $N_Y(x_2)$ is nonadjacent to any other vertices of $Y$. If
$|N_Y(x_2)| \geq m-2$ then form two sets $D_1$ and $D_2$ as follows.
By the maximality of $L$, one element in each pair $\{x_4,x_5\}$ and
$\{x_6,x_7\}$ is nonadjacent to all vertices of $N_Y(x_2)$. Call
these two vertices by $b$ and $c$. Therefore, there are at most four
edges connecting from $\{x_1,x_t\}$ to $\{x_3,b,c\}$ in $F$. Now,
define $D_1=\{x_1,x_t\} \cup \{$any $m-2$ vertices of $N_Y(x_2)\}$
and $D_2=\{x_3,b,c\} \cup$  $\{$any  $m-2$ other vertices of $Y\}$.
Thus, the set $D_1 \cup D_2$ induces $K_{5,6}$ minus four edges in $\overline{F}$, and so $\overline{F} \supseteq J_{10}$.\\

\noindent If $|N_Y(x_2)|= m-3$ then By the maximality of $L$, one
vertex in $\{x_4,x_5\}$ is nonadjacent to all vertices of
$N_Y(x_2)$. Call this vertex by $b$. Therefore, there are at most
four edges connecting from $\{x_1,x_2,x_t\}$ to $\{x_3,b\}$ in $F$.
Now, take $D_1=\{x_1,x_2,x_t\} \cup N_Y(x_2)$, and $D_2 =\{x_3,b\}
\cup \{$any $m-1$ other vertices of $Y\}$. Then, $D_1 \cup D_2$
contains $K_{5,6}$ minus at most four edges
in $\overline{F}$. Therefore,  $\overline{F} \supseteq J_{10}$.\\

\noindent if $|N_Y(x_2)|= m-4$ then take
$D_1=\{x_1,x_2,x_{t-1},x_t\} \cup N_Y(x_2) \cup N_Y(x_{t-1})$, and
$D_2 = \{x_3\} \cup \{$all the remaining vertices of $Y\}$. Then,
$D_1 \cup D_2$ contains $K_{5,6}$ minus possibly two edges
$(x_3,x_{t-1})$ and $(x_3,x_t)$
in $\overline{F}$. Therefore,  $\overline{F} \supseteq J_{10}$.\\

\noindent if $|N_Y(x_2)|= m-5$ then take
$D_1=\{x_1,x_2,b,x_{t-1},x_t\}$ where $b$ is either $x_3$ or $x_4$
whose the smallest number of neighbors in $Y$, and $D_2=Y$. Then,
$D_1 \cup D_2$ contains $K_{5,6}$ minus at most three edges in
$\overline{F}$. Therefore,  $\overline{F} \supseteq J_{10}$. \qed

\vspace{0.5cm}
\noindent {\bf The proof of Theorem 2}.\\
For $n=4$ and $k=2$, consider  graph $G = K_1\cup K_7$. Clearly  $G$
contains no $2P_n$ and
 $\overline{G}$ contains no $J_4$. Hence $ R(2P_4,J_4) \geq 9$. To prove the upper bound,
 consider now graph $F$ of order 9 containing no $2P_4$.
Take a longest path in $F$ and call it $L$. Let $L$ be $x_1,x_2,
\cdots, x_k$. Clearly, $k \leq 7$, since $F \not\supseteq 2P_4$. If
$A= V(F)\backslash V(L)$,
 then $|A| \geq 2$. Any vertex of $A$ is nonadjacent to $x_1$ and $x_k$.
  Thus, the number vertices in $A$ must be exactly 2 and so $k=7$,
  since otherwise $A$ together with $\{x_1,x_k\}$ will form a
  $K_{2,3}=J_4$ in $\overline{F}$. Let $A=\{y,z\}$, and consider the following two cases:\\

 \noindent {\bf Case 1.} Vertices $y$ and $z$ has a common neighbor in $L$. \\
 \noindent Let $x_i$ be the common neighbor of $y$ and $z$ in $L$, for some $i \in \{2,3, \cdots, 6\}$. Then,  $y,z$ are nonadjacent to  $x_{i-1}$ and $x_{i+1}$, since otherwise the maximality of $L$ will suffer. At least one of the last two vertices must differ with $x_1$ and $x_7$, call it $w$. So, we  have a $J_4$ in $\overline{F}$ formed by $\{x_1,x_7,y,z,w\}$.\\

\noindent {\bf Case 2.}  Vertices $y$ and $z$ has no common neighbor in $L$. \\
\noindent If there exists a vertex $x_i$, $2 \leq i \leq 6$, is nonadjacent to $y$ and $z$, then $\{x_i,x_1,x_7,y,z\}$ forms a $J_4$ in $\overline{F}$. Thus, every $x_i$ is adjacent to at least one of $\{y,z\}$. Now, since $y$ and $z$ has no common neighbor in $L$, without loss of generality we can assume that $x_2y \in E(F)$, and so $x_2z \notin E(F),  x_3y\notin E(F), x_3z\in E(F), x_4z\notin E(F), x_4y\in E(F), x_5y\notin E(F)$ and $x_5z\in E(F)$. Therefore, the path $x_1,x_2,y,x_4,x_3,z,x_5,x_6,x_7$ is Hamiltonian, which contradicts the maximality of path $L$ in $F$. \\

\noindent Now, let $n \geq 5$. Consider graph $G = K_1\cup
K_{kn-1}$. Clearly  $G$ contains no $kP_n$ and $\overline{G}$
contains no $J_4$. Hence $R(kP_n,J_4)\geq kn-1+1+1 =kn+1$. For the
upper bound, let $F$ be a graph of order $kn+1$ such that
$\overline{F}$ does not contain $J_4$.  By induction on $k$, we will
show that ${F}$ contains $kP_n$. By Theorem $A$ gives a verification
of the result for $k=1$. Assume the theorem is true for any $s \leq
k-1$, namely $R(sP_n, J_4)=sn+1$, for $n \geq 5$. Now consider graph
$F$ of $kn+1$ vertices such that $\overline{F} \not\supseteq J_4$.
By the induction hypothesis, $F$ will contain $(k-1)P_n$. Let
$Y=V(F)\backslash V((k-1)P_n)$. Then, $|Y|=n+1=R(P_n,J_4)$ and hence
$F[Y]$ contains a $P_n$. In total, $F$ will contain $kP_n$. \qed

\vspace{0.5cm}
\noindent {\bf The proof of Theorem 3.}\\
\noindent  Since graph $G = K_{m-1}\cup K_{kn-1}$ contains no $kP_n$
and $\overline{G}$ contains no $J_{2m}$, then $R(kP_n,J_{2m})\geq
kn+m-1$. For proving the upper bound, let $F$ be a graph of order
$kn+m-1$ such that  $\overline{F}$ contains no a $J_{2m}$.  We will
show that $F$ contains $k P_n$. We use an induction on $k$. For
$k=1$ it is true from Theorem $A$. Now, let assume that the theorem
is true for all $s \leq k-1$. Take any graph $F$ of $kn+m-1$
vertices such that its complement contains no $J_{2m}$. By the
hypothesis, $F$ must contain $(k-1)$ disjoint copies of $P_n$.
Remove these copies from $F$, then the remaining vertices will
induce another $P_n$ in $F$ since $\overline{F} \not\supseteq
J_{2m}$. Therefore $F \supseteq kP_n$. The proof is complete. \qed


\begin{thebibliography}{1}



\bibitem{BBS:05} E. T. Baskoro and Surahmat, The
Ramsey numbers of path with respect to wheels, \emph{Discrete
Math.}, \textbf{294} (2005), 275-277.


\bibitem{Bur:75} S. A. Burr, P. Erd\"{o}s and J. H. Spencer,  Ramsey theorem for
multiple copies of graphs, {\em Trans. Amer. Math. Soc.},
\textbf{209} (1975), 87-89.

\bibitem{BRR:88} S. A. Burr, On Ramsey numbers for large disjoint unions of graphs,
\emph{Discrete Math.}, \textbf{70} (1988), 277-293.

\bibitem{CZZ:04} Y. J. Chen, Y. Q. Zhang and K.M. Zhang, The Ramsey numbers of stars versus wheels,
 {\em European J. Combin.}, \textbf{25} (2004), 1067-1075.

 \bibitem{Catal:72}  V. Chv\'atal and F. Harary, Generalized Ramsey theory for
graphs, III: small off-diagonal numbers, \textit{Pac. J. Math.},
\textbf{41} (1972), 335-345.


\bibitem{HasB:04} Hasmawati, E. T. Baskoro and H. Assiyatun, Star-wheel Ramsey numbers,
\textit{J. Combin. Math. Conbin. Comput.}, \textbf{55} (2005),
123-128.


\bibitem{Rad:04}  S. P. Radziszowski, Small Ramsey numbers, {\em Electron.
J. Combin.}, July (2004) \#DS1.9, http://www.combinatorics.org/

\bibitem{SB:01}  Surahmat and E. T. Baskoro, On the Ramsey number of a
path or a star versus $W_4$ or $W_5$, \textit{Proceedings of the
12-th Australasian Workshop on Combinatorial Algorithms}, Bandung,
Indonesia, July 14-17 (2001), 165-170.

\bibitem{ST:06} Surahmat
and Ioan Tomescu, On path-Jahangir Ramsey numbers, {\it Preprint}
(2006).

\end{thebibliography}
\end{document}